# SELFEXTENSIONS OF MODULES OVER GROUP ALGEBRAS

KARIN ERDMANN, VIKTÓRIA KLÁSZ, AND RENÉ MARCZINZIK,
WITH AN APPENDIX BY BERNHARD BÖHMLER AND RENÉ MARCZINZIK

Abstract. Let $KG$ be a group algebra with $G$ a finite group and $K$ a field and $M$ an indecomposable $KG$-module. We pose the question, whether $\mathrm{Ext}^1_{KG}(M,M) \neq 0$ implies that $\mathrm{Ext}^i_{KG}(M,M) \neq 0$ for all $i \geq 1$. We give a positive answer in several important special cases such as for periodic groups and give a positive answer also for all Nakayama algebras, which allows us to improve a classical result of Gustafson. We then specialise the question to the case where the module $M$ is simple, where we obtain a positive answer also for all tame blocks of group algebras. For simple modules $M$, the appendix provides a Magma program that gives strong evidence for a positive answer to this question for groups of small order.

## Introduction

The cohomology theory of groups is a classical topic in homological algebra that has many applications to areas such as representation theory, algebraic topology and algebraic number theory, see for example [AM, Ben, Br]. One of the main goals of group cohomology is to study the behaviour of the Ext-groups $\mathrm{Ext}^i_{KG}(M,N)$ for $KG$-modules $M$ and $N$ and related cohomology rings when $G$ is a group and $K$ is a commutative ring. In this article we assume that $K$ is always a field, all groups are finite and all algebras are finite dimensional over a field $K$. We study the following question:

**Question M 0.1.** *Let $M$ be an indecomposable module over the group algebra $KG$ with $\mathrm{Ext}^1_{KG}(M,M) \neq 0$. Does this imply $\mathrm{Ext}^i_{KG}(M,M) \neq 0$ for all $i \geq 1$?*

We say that Question M has a positive answer for $KG$, if $\mathrm{Ext}^1_{KG}(M,M) \neq 0$ implies $\mathrm{Ext}^i_{KG}(M,M) \neq 0$ for all $i \geq 1$ for all indecomposable $KG$-modules $M$. We present the following evidence for Question M in this article:

**Theorem 0.2.** *Let $G$ be a group with group algebra $KG$ over a field $K$ such that one of the following conditions is satisfied:*

*(1) $G$ is a $p$-group, where $p$ denotes the characteristic of $K$.*
*(2) $G$ is abelian.*
*(3) The group algebra $KG$ is a periodic algebra.*

*Then Question M has a positive answer for $KG$.*

Recall here that a group is called *periodic* if its integer group cohomology is periodic, that is there exists a positive integer $d$ such that $H^n(G,\mathbb{Z}) \cong H^{n+d}(G,\mathbb{Z})$ for all $n \geq 1$. Being periodic is equivalent to all abelian subgroups being cyclic by a classical result of Artin and Tate, see for example [CE, Chapter XII.11]. Periodic groups also play an important role in algebraic topology, see for example [S]. It was shown by Erdmann and Skowroński (see for example [ES, Chapter 4] for a survey) that periodic groups are exactly the finite groups such that the group algebra $KG$ is a periodic algebra for all fields $K$, which means that all indecomposable non-projective $KG$-modules satisfy $\Omega^i(M) \cong M$ for some $i \geq 1$, see also [ES2] for equivalent characterisations of being periodic for group algebras. We note that the case of periodic $KG$ includes all cases where $KG$ is a representation-finite algebra. In fact, we will deduce Question M for representation-finite group algebras from the following more general result, using that all representation-finite group algebras are Brauer tree algebras and thus stable equivalent to a symmetric Nakayama algebra by classical results in [G] and [GR].

**Theorem 0.3.** *Let $A$ be a Nakayama algebra and $M$ an indecomposable $A$-module. Then* $\mathrm{Ext}^1_A(M,M) \neq 0$ *implies* $\mathrm{Ext}^i_A(M,M) \neq 0$ *for all* $i \geq 1$.

We apply the previous theorem also to obtain an improved version of an old result of Gustafson in [Gus] about bounds of the Loewy length for Nakayama algebras of finite global dimension. We cannot test Question M for all indecomposable modules over a given group algebra $KG$ as those algebras are in general wild. We can, however, test it for the finitely many simple modules. This leads us to the following conjecture:

**Conjecture S 0.4.** *Let $S$ be a simple $KG$ module with* $\mathrm{Ext}^1_{KG}(S,S) \neq 0$. *Then* $\mathrm{Ext}^i_{KG}(S,S) \neq 0$ *for all* $i \geq 1$.

For simple modules, the condition $\mathrm{Ext}^1_{KG}(S,S) \neq 0$ has a very visual interpretation: It means that there are loops in the quiver of the algebra $KG$ at the vertex corresponding to $S$. The strong no loop conjecture for finite dimensional algebras states that for an algebra $A$ and a simple $A$-module $S$, $\mathrm{Ext}^1_A(S,S) \neq 0$ implies that the projective dimension of $S$ is infinite. This was proved recently for finite dimensional algebras over algebraically closed fields in [ILP]. The open extreme no loop conjecture, see [LM], states that $\mathrm{Ext}^1_A(S,S) \neq 0$ even implies $\mathrm{Ext}^i_A(S,S) \neq 0$ for infinitely many $i$ and our Conjecture S can be seen as an even stronger conjecture for group algebras. Another related conjecture, see for example [GKM], is that for all simple modules over the group algebra $KS_n$ of the symmetric group $S_n$ for a field $K$ not of characteristic two, we have $\mathrm{Ext}^1_{KS_n}(S,S) = 0$ for all simple modules. The truth of our Conjecture S, would show that it is enough to show $\mathrm{Ext}^i_{KS_n}(S,S) = 0$ for some $i > 0$ to conclude this conjecture for the group algebra of the symmetric group. We partially tested Conjecture S for all group algebras $KG$ where $G$ has order at most 1000 over a splitting field using Magma where we showed that $\mathrm{Ext}^1_{KG}(S,S) \neq 0$ implies $\mathrm{Ext}^i_{KG}(S,S) \neq 0$ for all $i$ with $1 \leq i \leq 3$. The Magma program for this is presented in the appendix of this article. We also present the following theoretical evidence for conjecture $S$ in this article:

**Theorem 0.5.** *Let $A$ be a tame block of a group algebra. Then conjecture $S$ is true for every simple module in the block $A$.*

Here we include also the representation-finite algebras in the class of tame algebras by definition. Furthermore, we show that conjecture $S$ is true for all group algebras $KG$ over a field of characteristic $p$ such that $G$ has a normal $p$-complement, a result that helped a lot in making a quick program with Magma to test conjecture $S$.

## Acknowledgments

We thank Gjergji Zaimi for allowing us to include his quick proof of (1) in Proposition 2.3. Rene Marczinzik is grateful to Apolonia Gottwald for useful conversations and explaining an alternative proof of an inequality for $\mathrm{Ext}^1$ in Nakayama algebras to him. We also thank Jeremy Rickard for bringing some interesting examples to our attention. The results of this article are motivated by experiments with the GAP-package QPA, see [QPA], and Magma, see [BCP]. We thank Johannes Schmitt for useful comments on the Magma code. This work is partly supported by the Deutsche Forschungsgemeinschaft (DFG, German Research Foundation) under Germany's Excellence Strategy - GZ 2047/1, Projekt-ID 390685813.

## 1. Selfextensions for modules over Nakayama and related algebras

We assume that all algebras are finite dimensional $K$-algebras over a field $K$ and modules are right modules unless stated otherwise. Let $J$ denote the Jacobson radical of an algebra and $D = \mathrm{Hom}_K(-,K)$ the natural duality of a finite dimensional algebra $A$. A finite dimensional algebra $A$ is a Nakayama algebra if every indecomposable projective left and right module is uniserial, meaning it has a unique composition series. The homological questions we study in this article are invariant under Morita-equivalence and we are mainly interested in the questions over splitting fields for the algebras. Thus, for the case of Nakayama algebras we can and will assume that they are given by a quiver and admissible relations. The quiver of a connected Nakayama algebra is then either a linearly oriented line:

$$\bullet^0 \longrightarrow \bullet^1 \longrightarrow \cdots \longrightarrow \bullet^{n-2} \longrightarrow \bullet^{n-1}$$

or a linearly oriented cycle:

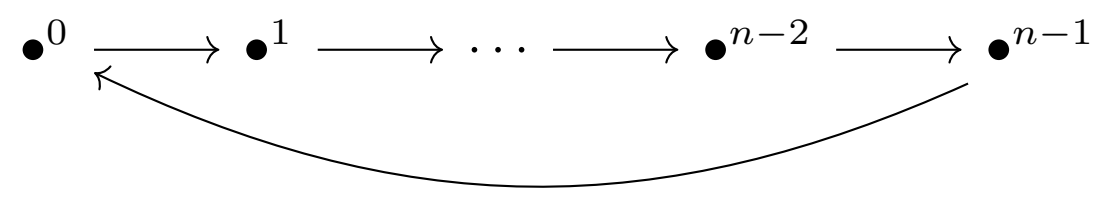

Recall that a module $M$ over an algebra $A$ is called rigid if $\mathrm{Ext}_A^1(M,M)=0$ and we call it non-rigid if $\mathrm{Ext}_A^1(M,M)\neq 0$. In this section, we give an elementary translation when an indecomposable module over a Nakayama algebra is non-rigid and use that to prove our main results. We assume that the reader is familiar with the basics of representation theory of algebras and homological algebra. We refer to the books [ARS] and [AnFul] for chapters on Nakayama algebras, and to [Mar] on the calculation of projective resolutions for Nakayama algebras. For a given Nakayama algebra $A$, we denote by $c_i$ the Loewy length of the indecomposable projective module $e_iA$ when fixing a complete set of primitive orthogonal idempotents $e_i$. The sequence $[c_0, c_1, ..., c_{n-1}]$ is called the Kupisch series of a Nakayama algebra $A$ and determines the algebra uniquely (assuming it is given by a connected quiver algebra over a fixed field). Note that Nakayama algebras are uniquely characterised by the condition that every indecomposable module is uniserial, and thus every indecomposable module over a Nakayama algebra $A$ can be written in the form $e_iA/e_iJ^k$ for some $k$ with $1\le k\le c_i$. As usual, we assume that the $c_i$ are defined for all $i\in\mathbb{Z}$ by setting $c_i = c_k$ if $i=k$ modulo $n$, where $n$ denotes the number of simple modules of the algebra. We will often use the following lemma to calculate Hom-spaces for Nakayama algebras.

**Lemma 1.1.** *Let $A$ be a Nakayama algebra. Then* $\mathrm{Hom}_A(e_iA/e_iJ^k, e_jA/e_jJ^l)\cong (e_jJ^{\max(0,l-k)}/e_jJ^l)e_i$, *where we interpret $J^0=A$.*

*Proof.* See [Mar2, Lemma 1.5].

□

**Lemma 1.2.** *Let $A$ be a Nakayama algebra with $n$ simple modules and $M=e_iA/e_iJ^k$ and $N=e_sA/e_sJ^t$ indecomposable, non-projective $A$-modules.*

*(1) If $t\ge k$, we have* $\mathrm{Ext}_A^1(N,M)\cong\mathrm{Hom}_A(\Omega^1(N),M)$.

*(2) $M$ is non-rigid if and only if $n\le k\le c_i-n$.*

*Proof.* (1) Let $0\to\Omega^1(N)\to P\to N\to 0$ be the short exact sequence such that $P\to N\to 0$ is a projective cover of $N$. We apply the functor $\mathrm{Hom}_A(-,M)$ to this short exact sequence to obtain:

$$0\to\mathrm{Hom}_A(N,M)\xrightarrow{f}\mathrm{Hom}_A(P,M)\xrightarrow{g}\mathrm{Hom}_A(\Omega^1(N),M)\xrightarrow{h}\mathrm{Ext}_A^1(N,M)\to 0.$$

Note that $\mathrm{Hom}_A(N,M)=\mathrm{Hom}_A(e_sA/e_sJ^t, e_iA/e_iJ^k)$ is given by left multiplication maps of the form $l_z$ with $z\in(e_iA/e_iJ^k)e_s$ such that $zJ^t=0$. But since we assume that $t\ge k$, the condition $zJ^t=0$ is automatic for $z\in(e_iA/e_iJ^k)e_s$.

Thus, $\mathrm{Hom}_A(N,M)\cong(e_iA/e_iJ^k)e_s$ and $\mathrm{Hom}_A(P,M)=\mathrm{Hom}_A(e_sA,e_iA/e_iJ^k)=(e_iA/e_iJ^k)e_s$ have the same length. Since $f$ is an injective map between modules of the same length, it is an isomorphism. The above exact sequence thus gives us $\mathrm{Ker}(g)=\mathrm{Im}(f)=\mathrm{Hom}_A(P,M)$ and thus $g=0$. This shows that $\mathrm{Ker}(h)=\mathrm{Im}(g)=0$ and thus $h$ is injective, but by the above exact sequence it is also surjective. Thus $h$ is an isomorphism and we obtain $\mathrm{Ext}_A^1(N,M)\cong\mathrm{Hom}_A(\Omega^1(N),M)$.

(2) By (1) applied with $N=M$, $M$ is non-rigid if and only if $\mathrm{Hom}_A(\Omega^1(M),M)\neq 0$. Note that $\Omega^1(M)\cong e_{i+k}A/e_{i+k}J^{c_i-k}$ and thus $\mathrm{Hom}_A(\Omega^1(M),M)\cong\mathrm{Hom}_A(e_{i+k}A/e_{i+k}J^{c_i-k},e_iA/e_iJ^k)$. Now note that by Lemma 2.1: $\mathrm{Hom}_A(e_{i+k}A/e_{i+k}J^{c_i-k},e_iA/e_iJ^k)\cong(e_iJ^{\max(0,2k-c_i)}/e_iJ^k)e_{i+k}$, where we interpret $J^0=A$. Now $(e_iJ^{\max(0,2k-c_i)}/e_iJ^k)e_{i+k}\neq 0$ is equivalent to $k\ge n$ and $\max(0,2k-c_i)\le s+(t-1)n=k-n$ when $k=s+tn$ for some $s$ with $0\le s\le n-1$ and $t\ge 1$. Now note that the conditions $n\le k$ and $\max(0,2k-c_i)\le k-n$ are equivalent to $n\le k\le c_i-n$.

□

**Proposition 1.3.** *Let $A$ be a Nakayama algebra and $M$ an indecomposable, non-projective $A$-module.*

*(1) If $M$ is non-rigid, then $\Omega^1(M)$ is also non-rigid.*

*(2) If $M$ is non-rigid, it has infinite projective dimension and infinite injective dimension.*

*Proof.* Assume $A$ has $n$ simple modules.

(1) Assume $M = e_iA/e_iJ^k$ is non-rigid, which by Lemma 2.2 (2) is equivalent to $n \leq k \leq c_i - n$. By Lemma 2.2 (2), $\Omega^1(M) = e_{i+k}A/e_{i+k}J^{c_i-k}$ is non-rigid if and only if $n \leq c_i - k \leq c_{i+k} - n$. We thus have to show that $n \leq k \leq c_i - n$ implies $n \leq c_i - k \leq c_{i+k} - n$. The inequality $n \leq c_i - k$ follows immediately from the assumption that $M$ is non-rigid. The second inequality $c_i - k \leq c_{i+k} - n$ follows from observing that

$$c_{i+k} - c_i = c_{i+k} - c_{i+n} = \sum_{j=n+1}^{k} (c_{i+j} - c_{i+j-1}) \geq \sum_{j=n+1}^{k} (-1) = n - k$$

(here we used that $c_l - c_{l-1} \geq -1$ for a general Kupisch series of a Nakayama algebra). After rearranging the terms we get $c_i - k \leq c_{i+k} - n$ as desired.

(2) Assume $M$ is non-rigid. By (1) of this theorem also $\Omega^1(M)$ is non-rigid and thus by induction $\Omega^k(M)$ for all $k \geq 1$ is non-rigid and thus, especially non-projective. This gives that $M$ has infinite projective dimension. To see that $M$ also has infinite injective dimension, note that $\mathrm{Ext}^1_{A^{op}}(D(M), D(M)) \cong \mathrm{Ext}^1_A(M, M) \neq 0$ and thus $D(M)$ has infinite projective dimension, which is equivalent to $M$ having infinite injective dimension by applying the duality $D$.

□

**Corollary 1.4.** *Let $A$ be a Nakayama algebra with finite global dimension. Then $\mathrm{Ext}^1_A(M, M) = 0$ for all indecomposable $A$-modules $M$.*

We now prove our first main result.

**Theorem 1.5.** *Let $A$ be a Nakayama algebra with indecomposable $A$-modules $M, N$. Then the following hold:*

*(1)*

$$\dim(\mathrm{Hom}_A(\Omega^1(M), N)) \geq \min(\dim(\mathrm{Hom}_A(\Omega^1(M), M)), \dim(\mathrm{Hom}_A(\Omega^1(N), N)))$$

*(2)*

$$\dim(\mathrm{Ext}^1_A(M, N)) \geq \min(\dim(\mathrm{Ext}^1_A(M, M)), \dim(\mathrm{Ext}^1_A(N, N)))$$

*Proof.* (1) As discussed at the beginning of this chapter, we can assume that the Nakayama algebra $A$ is given as the path algebra of a quiver $Q$ modulo an admissible ideal. If $A$ has $n$ simple modules then $Q$ is either an oriented cycle or an oriented line with $n$ vertices, denoted by $0, 1, \ldots, n-1$. In the following, when we refer to an integer $i$ as a vertex of $Q$ then we mean the unique vertex $j \in \{0, 1, \ldots, n-1\}$ which is congruent to $i$ modulo $n$.

Let $M = e_iA/e_iJ^k$ and $N = e_sA/e_sJ^t$. Then $\Omega^1(M) = e_{i+k}A/e_{i+k}J^{(c_i-k)}$ and $\Omega^1(N) = e_{s+t}A/e_{s+t}J^{(c_s-t)}$, so Lemma 2.1 implies that $d_1 := \dim(\mathrm{Hom}_A(\Omega^1(N), N)) = \dim((e_sJ^{\max(0,2t-c_s)}/e_sJ^t)e_{s+t})$ is equal to the number of paths $p$ from vertex $s$ to vertex $s+t$ in the quiver $Q$ such that $\max(0, 2t - c_s) \leq \mathit{length}(p) < t$. There are exactly $\lfloor \frac{t}{n} \rfloor$ paths from $s$ to $s + t$ of length at most $t - 1$. If $Q$ is an oriented line, it is clear that there are 0 such paths and $0 = \lfloor \frac{t}{n} \rfloor$ in this case. If $Q$ is an oriented cycle and $t = n \cdot m + d$ with $m, d \in \mathbb{N}_0$ and $0 \leq d < n$ then there are $m = \lfloor \frac{t}{n} \rfloor$ such paths, depending on how many times we go around; namely, one path of length $d$, one of length $n + d$, one of length $2n + d, \ldots$ and one of length $(m-1) \cdot n + d$. Thus, $d_1 \leq \lfloor \frac{t}{n} \rfloor$.

Similarly, using Lemma 2.1 we obtain that $d_2 := \dim(\mathrm{Hom}_A(\Omega^1(M), N)) = \dim((e_sJ^{\max(0,t-c_i+k)}/e_sJ^t)e_{i+k})$ is equal to the number of paths $p$ from $s$ to $i + k$ with $\max(0, t - c_i + k) \leq \mathit{length}(p) < t$. If $\max(0, t - c_i + k) = 0$ then $d_2 \geq \lfloor \frac{t}{n} \rfloor$. In this case, we obtain that $d_2 \geq \lfloor \frac{t}{n} \rfloor \geq d_1$, which implies the desired inequality.

It remains to prove the statement in the case when $t > c_i - k$. In this case, there are at least $\lfloor \frac{t-(t-c_i+k)}{n} \rfloor = \lfloor \frac{c_i-k}{n} \rfloor$ paths $p$ from $s$ to $i + k$ with $\max(0, t - c_i + k) \leq \mathit{length}(p) < t$ in $Q$. So, $d_2 \geq \lfloor \frac{c_i-k}{n} \rfloor$. Moreover, there are always at most $\lfloor \frac{c_i-k}{n} \rfloor$ many paths $p$ from vertex $i$ to vertex $i + k$ in $Q$ with $\max(0, 2k - c_i) \leq \mathit{length}(p) < k$. Indeed, if $2k - c_i \geq 0$ then there are exactly $\lfloor \frac{k-(2k-c_i)}{n} \rfloor = \lfloor \frac{c_i-k}{n} \rfloor$ many such paths; and if $k \leq c_i - k$ then $\lfloor \frac{k}{n} \rfloor$ many, and in this case we have $\lfloor \frac{k}{n} \rfloor \leq \lfloor \frac{c_i-k}{n} \rfloor$. Thus, $d_3 := \dim(\mathrm{Hom}_A(\Omega^1(M), M)) = \dim((e_iJ^{\max(0,2k-c_i)}/e_iJ^k)e_{i+k}) \leq \lfloor \frac{c_i-k}{n} \rfloor$. In conclusion, $d_2 \geq d_3$, which proves the claim.

(2) Assume first that $\dim M \geq \dim N$. Then we know that by Lemma 2.2 (1): $\mathrm{Ext}_A^1(M,N) \cong \mathrm{Hom}_A(\Omega^1(M),N)$, $\mathrm{Ext}_A^1(M,M) \cong \mathrm{Hom}_A(\Omega^1(M),M)$ and $\mathrm{Ext}_A^1(N,N) \cong \mathrm{Hom}_A(\Omega^1(N),N)$. Now by (1) of this theorem we also know that

$$\dim \mathrm{Hom}_A(\Omega^1(M),N) \geq \min(\dim \mathrm{Hom}_A(\Omega^1(M),M), \dim \mathrm{Hom}_A(\Omega^1(N),N)),$$

which proves the claim when $\dim M \geq \dim N$.
Now assume that $\dim N \geq \dim M$. Then $\dim \mathrm{Ext}_A^1(M,N) = \dim \mathrm{Ext}_{A^{op}}^1(D(N),D(M))$ and since $\dim D(N) \geq \dim D(M)$, we can apply the first case and obtain $\dim \mathrm{Ext}_{A^{op}}^1(D(N),D(M)) \geq \min(\dim \mathrm{Ext}_{A^{op}}^1(D(M),D(M)), \dim \mathrm{Ext}_{A^{op}}^1(D(N),D(N)))$, which proves the claim also in this case since $\dim \mathrm{Ext}_{A^{op}}^1(D(X),D(X)) = \dim \mathrm{Ext}_A^1(X,X)$ for arbitrary $X$.

□

**Lemma 1.6.** *Let $A$ be a Nakayama algebra and $N$ an indecomposable $A$-module. For all $i \geq 1$:*

$$\dim(\mathrm{Ext}_A^i(N,N)) \geq \min(\dim(\mathrm{Ext}_A^1(\Omega^{i-1}(N),\Omega^{i-1}(N))), \dim(\mathrm{Ext}_A^1(N,N))).$$

*Proof.* Just set $M = \Omega^{i-1}(N)$ in Theorem 2.5 and use that $\mathrm{Ext}_A^i(X,Y) \cong \mathrm{Ext}_A^1(\Omega^{i-1}(X),Y)$ in general. □

**Theorem 1.7.** *Let $A$ be a Nakayama algebra and $N$ an indecomposable $A$-module. If $\mathrm{Ext}_A^1(N,N) \neq 0$, we have $\mathrm{Ext}_A^i(N,N) \neq 0$ for all $i > 0$.*

*Proof.* Assume $\mathrm{Ext}_A^1(N,N) \neq 0$. By Proposition 2.3 (1), we have that with $N$ also all $\Omega^i(N)$ are non-rigid for $i \geq 1$. Then by Lemma 2.6, we have

$$\dim(\mathrm{Ext}_A^i(N,N)) \geq \min(\dim(\mathrm{Ext}_A^1(\Omega^{i-1}(N),\Omega^{i-1}(N))), \dim(\mathrm{Ext}_A^1(N,N))) > 0.$$

□

We apply our main result to give an improved result on the Loewy length of Nakayama algebras of finite global dimension that was first obtained in [Gus].

**Proposition 1.8.** *Let $A$ be a Nakayama algebra with Loewy length $L(A)$ and $n$ simple modules.*

*(1) Assume $L(A) \geq 2n$. Then there exists an indecomposable $A$-module $M$ with $\mathrm{Ext}_A^i(M,M) \neq 0$ for all $i > 0$.*
*(2) Assume $A$ has finite global dimension. Then $L(A) \leq 2n-1$.*

*Proof.* (1) Assume $A$ has Loewy length at least $2n$ and let $P = e_iA$ be an indecomposable projective module with Loewy length $c_i \geq 2n$. Let $M := e_iA/e_iJ^k$ with $k = n$. Then we clearly have $n \leq k \leq c_i - n$ and thus by Lemma 2.2 (2) $M$ is non-rigid and thus $\mathrm{Ext}_A^i(M,M) \neq 0$ for all $i > 0$ by Theorem 2.6.
(2) Assume to the contrary that $A$ has finite global dimension but $L(A) \geq 2n$. Then by (1), $A$ has an indecomposable non-rigid module $M$, which contradicts Corollary 2.4. Thus we need to have $L(A) \leq 2n-1$.

□

By the previous result, we know that a Nakayama algebra $A$ with $n$ simple modules and $\mathrm{Ext}_A^1(M,M) = 0$ for all indecomposable $A$-modules $M$ has Loewy length at most $2n-1$. This is the perfect bound for Nakayama algebras as the Nakayama algebra with Kupisch series $[n, 2n-1, 2n-2, ..., n+1]$ has finite global dimension and Loewy length $2n-1$. Indeed, this Nakayama algebra is then given by $A = kQ/I$, where $Q$ is the quiver

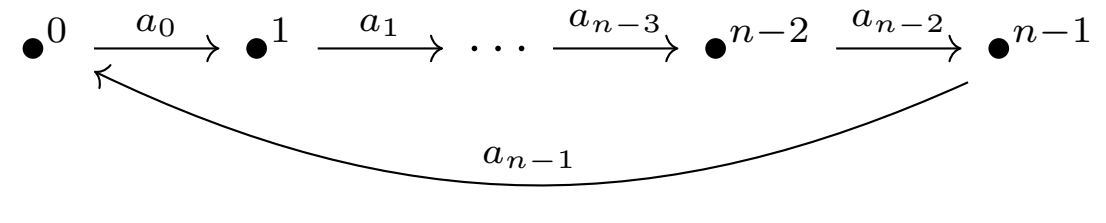

and $I = \langle a_0 a_1 \dots a_{n-1} \rangle$. The indecomposable projective modules in this case look like

$$P(0) = \begin{smallmatrix} 0 \\ 1 \\ \vdots \\ n-1 \end{smallmatrix}, \; P(1) = \begin{smallmatrix} 1 \\ 2 \\ \vdots \\ n-1 \\ 0 \\ \vdots \\ n-1 \end{smallmatrix}, \; P(2) = \begin{smallmatrix} 2 \\ 3 \\ \vdots \\ n-1 \\ 0 \\ \vdots \\ n-1 \end{smallmatrix}, \; \ldots, \; P(n-1) = \begin{smallmatrix} n-1 \\ 0 \\ \vdots \\ n-1 \end{smallmatrix}$$

and have Loewy length equal to $n, 2n-1, 2n-2, \ldots, n+1$, respectively. Because the Loewy length of $A$, is equal to the maximum of the Loewy lengths of the projectives, $L(A) = 2n-1$. In order to see that the global dimension of $A$ is finite, we only need to show that the simple modules $S_0, \ldots, S_{n-1}$ have finite projective dimension. The corresponding minimal projective resolutions look like

$$0 \longrightarrow P(i+1) \longrightarrow P(i) \longrightarrow S_i \longrightarrow 0$$

for every $i \neq 0$ where we set $P(n) = P(0)$, and

$$0 \longrightarrow P(0) \longrightarrow P(1) \longrightarrow P(0) \longrightarrow S_0 \longrightarrow 0.$$

So $\mathrm{gldim}(A) = 2 < \infty$.

We remark that it is an open problem whether an algebra $A$ with $\mathrm{Ext}^1_A(M, M) = 0$ for all indecomposable modules $M$ is always representation-finite. This is true over algebraically closed fields, see for example the end of the article [Sko]. The next example shows that it is in general not true that an indecomposable module $M$ over a Nakayama algebra $A$ with $\mathrm{Ext}^i_A(M, M) \neq 0$ for some $i \geq 2$ has infinite projective dimension. We also give an example that shows that in general we have that $\mathrm{Ext}^n_A(M, M)$ is not isomorphic to $\mathrm{Hom}_A(\Omega^n(M), M)$ for Nakayama algebras for $n \geq 2$, despite this being true for $n = 1$ by Lemma 2.2 (1).

**Example 1.9.** Let $A$ be the Nakayama algebra given by a quiver and relations with Kupisch series $[2, 2, ..., 2, 3]$ with $n \geq 2$ entries where all but the last entry are equal to 2. This means that $A \cong kQ/I$, where $Q$ is the quiver and $I = \langle a_0a_1, \; a_1a_2, \; \ldots, \; a_{n-2}a_{n-1} \rangle$. The indecomposable projective modules in this case look like

$$P(0) = \begin{smallmatrix} 0 \\ 1 \end{smallmatrix}, \; P(1) = \begin{smallmatrix} 1 \\ 2 \end{smallmatrix}, \; \ldots, \; P(n-2) = \begin{smallmatrix} n-2 \\ n-1 \end{smallmatrix}, \; P(n-1) = \begin{smallmatrix} n-1 \\ 0 \\ 1 \end{smallmatrix}$$

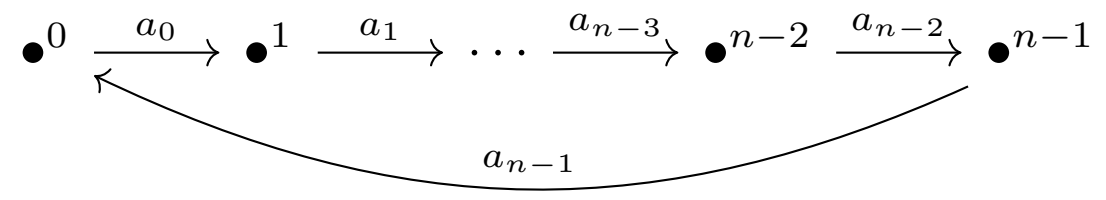

$A$ has finite global dimension equal to $n$ and there is a unique simple module $S_0$ with projective dimension equal to $n$ but $\mathrm{Ext}^n_A(S_0, S_0) \neq 0$. There also exists a unique indecomposable module $M = \begin{smallmatrix} n-1 \\ 0 \end{smallmatrix}$ with dimension vector $[1, 0, ..., 0, 1]$ (with exactly two non-zero entries) such that $\mathrm{Ext}^n_A(M, M)$ is zero but $\mathrm{Hom}_A(\Omega^n(M), M)$ is one-dimensional.

The next example shows that Theorem 2.7 cannot be extended to general uniserial modules over general algebras.

**Example 1.10.** Let $Q$ be the Kroenecker quiver and $A = kQ$ the Kroenecker algebra. Then there exist two-dimensional uniserial $A$-modules which are non-rigid, see for example [ARS, Chapter VIII, Section 7]. But since $A$ is hereditary, it has finite global dimension and thus all those uniserial non-rigid modules have finite projective dimension.

## 2. Group algebras

We first use our result on Nakayama algebras from the previous section to show that Question $M$ has a positive answer for representation-finite blocks of group algebras. Recall that the Tate cohomology groups $\hat{\mathrm{Ext}}^i_A(M, N)$ for two modules $M$ and $N$ over a selfinjective algebra $A$ are defined as $\hat{\mathrm{Ext}}^i_A(M, N) := \underline{\mathrm{Hom}}_A(\Omega^i(M), N)$ for $i \in \mathbb{Z}$. We refer to [CTVZ, section 2.6] for more on the Tate cohomology and recall that for general selfinjective algebras $A$ and for $i \geq 1$, we have $\mathrm{Ext}_A^i(M,N) := \underline{\mathrm{Hom}_A}(\Omega^i(M),N)$, see for example [SkoYam, Chapter IV, Theorem 9.6].

In the rest of this article we will frequently use the following lemma:

**Lemma 2.1.** *Let $A$ be a finite dimensional algebra, $S$ a simple $A$-module and $M$ an $A$-module with minimal projective resolution*

$$\cdots \to P_n \to P_{n-1} \to \cdots \to P_0 \to M \to 0.$$

*Then* $\mathrm{Ext}_A^n(M,S) \cong \mathrm{Hom}_A(\Omega^n(M),S) \cong \mathrm{Hom}_A(P_n,S)$.

*Proof.* See [Ben, Corollary 2.5.4]. □

We now apply our results to selfinjective Nakayama algebras and Brauer tree algebras.

**Lemma 2.2.** *Let $A$ be a selfinjective Nakayama algebra and $N$ an indecomposable $A$-module. Then* $\dim(\hat{\mathrm{Ext}}_A^i(N,N)) \geq \dim(\mathrm{Ext}_A^1(N,N))$ *for all $i \in \mathbb{Z}$.*

*Proof.* First assume that $i$ is positive. By Lemma 2.6 we have

$$\dim(\mathrm{Ext}_A^i(N,N)) \geq \min(\dim(\mathrm{Ext}_A^1(\Omega^{i-1}(N),\Omega^{i-1}(N))),\dim(\mathrm{Ext}_A^1(N,N))).$$

Now use that $\Omega^i$ is an equivalence of the stable category of a selfinjective algebra to see that for all $i > 0$ we have $\mathrm{Ext}_A^1(\Omega^{i-1}(N),\Omega^{i-1}(N)) \cong \mathrm{Ext}_A^1(N,N)$ to obtain $\dim(\hat{\mathrm{Ext}}_A^i(N,N)) \geq \dim(\mathrm{Ext}_A^1(N,N))$ for all $i > 0$. Now assume that $i \leq 0$ and note that $N$ is a periodic module since $A$ is representation-finite selfinjective. Let $i-1 \equiv l \bmod m$, where $m$ is the period of $N$ for a positive $l$. Then $\hat{\mathrm{Ext}}_A^i(N,N) \cong \mathrm{Ext}_A^1(\Omega^{i-1}(N),N) \cong \mathrm{Ext}_A^1(\Omega^l(N),N)$ and $l$ is positive so the result follows from the first part of the proof. □

Recall that Brauer tree algebras are exactly the special biserial symmetric representation-finite algebras. For a more explicit definition in terms of quiver and relations, we refer for example to [SkoYam, Chapter IV, Section 4]. A well-known and important result for Brauer tree algebras is that they are exactly those symmetric algebras that are stable equivalent to symmetric Nakayama algebras, see for example [ARS, Chapter X.3].

**Theorem 2.3.** *Let $A$ be a Brauer tree algebra and $M$ an indecomposable $A$-module.*

*(1)* $\dim(\hat{\mathrm{Ext}}_A^i(M,M)) \geq \dim(\mathrm{Ext}_A^1(M,M))$ *for all $i \in \mathbb{Z}$.*
*(2) If* $\mathrm{Ext}_A^1(M,M) \neq 0$*, then* $\hat{\mathrm{Ext}}_A^i(M,M) \neq 0$ *for all $i > 0$.*

*Proof.* Let $F : \underline{\mathrm{mod}} - A \to \underline{\mathrm{mod}} - B$ be a stable equivalence where $B$ is a symmetric Nakayama algebra. By [ARS, Chapter X, Proposition 1.12. (b)], the stable equivalence $F$ commutes with $\Omega$, which means that $F\Omega_A \cong \Omega_B F$. Now let $M$ be an indecomposable $A$-module with $\mathrm{Ext}_A^1(M,M) \neq 0$. Note that since $\mathrm{Ext}^l(X,Y) \cong \underline{\mathrm{Hom}}(\Omega^l(X),Y)$ for general modules $X,Y$ in a symmetric algebra and since $F$ commutes with $\Omega$ we have that $F$ preserves Ext, meaning that $\hat{\mathrm{Ext}}_A^l(X,Y) \cong \hat{\mathrm{Ext}}_B^l(F(X),F(Y))$ in general. Thus we also have $\mathrm{Ext}_B^1(F(M),F(M)) \neq 0$ and by Theorem 2.7, we get $\mathrm{Ext}_B^l(F(M),F(M)) \neq 0$ for all $l > 0$ since $B$ is a symmetric Nakayama algebra. But then we also have $\mathrm{Ext}_A^l(M,M) \cong \mathrm{Ext}_B^l(F(M),F(M)) \neq 0$ for all $l > 0$ and $\dim(\hat{\mathrm{Ext}}_A^i(M,M)) \geq \dim(\mathrm{Ext}_A^1(M,M))$ for all $i \in \mathbb{Z}$ as this holds for symmetric Nakayama algebras by Lemma 3.2. □

**Corollary 2.4.** *Let $KG$ be a representation-finite group algebra with an indecomposable $KG$-module $M$. If* $\mathrm{Ext}_A^1(M,M) \neq 0$*, then* $\hat{\mathrm{Ext}}_A^i(M,M) \neq 0$ *for all $i > 0$.*

We now focus on periodic groups $G$. We need the following result:

**Theorem 2.5.** *Let $G$ be a periodic group. Then every block $B$ of $KG$ is either representation-finite or of quaternion type. If $B$ is of infinite representation type, then every module is at most 4-periodic.*

*Proof.* See for example [ES, Section 4]. □

**Proposition 2.6.** *Let $A$ be a symmetric algebra and $M$ a 4-periodic module. Then* $\mathrm{Ext}_A^1(M,M) \neq 0$ *implies* $\mathrm{Ext}_A^i(M,M) \neq 0$ *for all $i \geq 1$.*

*Proof.* Since $M$ is 4-periodic, we have that $\operatorname{Ext}_A^i(M,M) \cong \operatorname{Ext}_A^{i+4}(M,M)$ for all $i \geq 1$. Thus it is enough to show that $\operatorname{Ext}_A^i(M,M) \neq 0$ for $i = 2,3,4$ if $\operatorname{Ext}_A^1(M,M) \neq 0$. We have $\operatorname{Ext}_A^4(M,M) \cong \underline{\operatorname{Hom}}_A(\Omega^4(M),M) \cong \underline{\operatorname{Hom}}_A(M,M) \neq 0$, since $\Omega^4(M) \cong M$. Next we use the Auslander-Reiten formula $\operatorname{Ext}_A^1(M,N) \cong D\underline{\operatorname{Hom}}_A(N,\tau(M))$ and $\tau \cong \Omega^2$ for a symmetric algebra $A$. This gives us $\operatorname{Ext}_A^3(M,M) \cong \operatorname{Ext}_A^1(\Omega^2(M),M)) \cong D\underline{\operatorname{Hom}}_A(M,\tau(\Omega^2(M))) \cong D\underline{\operatorname{Hom}}_A(M,M) \neq 0$. Now note that in a symmetric algebra $\underline{\operatorname{Hom}}_A(M,N) \cong \underline{\operatorname{Hom}}_A(\Omega^i(M),\Omega^i(N))$ for all $i \in \mathbb{Z}$, since $\Omega$ is a stable equivalence. This gives us $\operatorname{Ext}_A^2(M,M) \cong \underline{\operatorname{Hom}}_A(\Omega^2(M),M) \cong \underline{\operatorname{Hom}}_A(M,\Omega^2(M)) \cong D\operatorname{Ext}_A^1(M,M) \neq 0$, where in the second isomorphism we applied $\Omega^2$ and used that $M$ is 4-periodic and in the last isomorphism we used the Auslander-Reiten formula again. □

**Corollary 2.7.** *Let $A$ be a symmetric algebra of quaternion type. Then Question $M$ has a positive answer for $A$. In particular, conjecture $M$ is true for periodic groups $G$.*

*Proof.* In a symmetric algebra $A$ of quaternion type we have that every indecomposable non-projective module is 4-periodic and thus by the previous proposition Question $M$ has a positive answer for $A$. The blocks of a group algebra $KG$ for periodic $G$ are either representation-finite or of quaternion type, see Theorem 3.5. In both cases we know that Question $M$ has a positive answer by Corollary 3.4 and the previous proposition, which proves the corollary.

□

The next example shows that for general symmetric algebras, Conjecture S can be false:

**Example 2.8.** Let $A = KQ/I$ be the algebra with the quiver $Q$ and relations $I$ as follows:

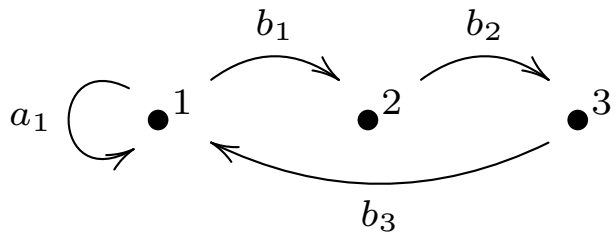

$I =< b_3b_1, a_1^2b_1, b_1b_2b_3 - a_1^2, b_3a_1^2, a_1^4, b_2b_3a_1b_1b_2 >$ .

$A$ is a representation-finite symmetric algebra and we have $\operatorname{Ext}_A^1(S,S) \neq 0$ but $\operatorname{Ext}_A^3(S,S) = 0$.

We now prove Question $M$ for local and commutative group algebras.

**Proposition 2.9.** *Let $A = KG$ be a group algebra with $G$ a $p$-group. Then $\operatorname{Ext}_A^i(M,M) \neq 0$ for all $i \geq 1$ and all non-projective modules $M$.*

*Proof.* Note that $A$ is a local algebra. Clearly it is enough to prove the statements for indecomposable $M$. Assume that $M$ is not projective. We have $\operatorname{Ext}_A^i(M,M) \cong \operatorname{Ext}_A^i(M \otimes_K K, M) \cong \operatorname{Ext}_A^i(K, \operatorname{Hom}_A(M,M))$ using [Ben, Proposition 3.1.8]. Furthermore, we have $\operatorname{Hom}_A(M,M) \cong M \otimes_K M^*$, which is projective if and only if $M$ is projective, see [Ben, Proposition 3.1.10]. Since $A$ is local and $\operatorname{Hom}_A(M,M)$ not projective, we have $\operatorname{Ext}_A^i(K, \operatorname{Hom}_A(M,M)) \neq 0$ as any non-projective module over a local algebra has infinite injective dimension. To see this, let $T$ be a general non-projective module over a general local finite dimensional algebra $A$, then $T$ has infinite projective dimension. We show that $T$ having finite projective dimension leads to a contradiction: Assume $T$ has finite projective dimension equal to $r$. Then $N := \Omega^{r-1}(T)$ has projective dimension equal to 1 and there exists a minimal projective presentation of the form

$$0 \to A^t \to A^s \to N \to 0.$$

This implies that $A^t$ is in the radical of $A^s$, which is a contradiction since both modules have the same Loewy length. Thus $T$ has infinite projective dimension. Thus if a general non-injective module $R$ is non-injective over a local algebra $A$, then $R$ has infinite injective dimension or else the non-projective module $D(R)$ would have finite non-zero projective dimension over the local algebra $A^{op}$. □

**Corollary 2.10.** *Let $G$ be abelian and $A = KG$ the group algebra of $G$. Let $M$ be an $A$-module with $\operatorname{Ext}_A^1(M,M) \neq 0$. Then $\operatorname{Ext}_A^i(M,M) \neq 0$ for all $i \geq 1$ and all non-projective modules $M$.*

*Proof.* Note that $G$ abelian implies that every block of $KG$ is also a group algebra of an abelian group and such an indecomposable block must either be semisimple or Morita equivalent to a local algebra. Thus the statement follows from the previous proposition. □

The most important case of Conjecture S is $S$ being the trivial module $K$. Therefore, we state the following special case of conjecture $S$:

**Conjecture K 2.11.** *Let $KG$ be a group algebra with $\mathrm{Ext}^1_{KG}(K,K) \neq 0$, then $\mathrm{Ext}^i_{KG}(K,K) \neq 0$ for all $i \geq 1$.*

**Proposition 2.12.** *Let $KG$ be a group algebra over an algebraically closed field $K$ with $\mathrm{Ext}^1_{KG}(K,K) \neq 0$ and assume conjecture K is true for $KG$. Then Question M has a positive answer for all indecomposable $KG$-modules $M$ such that $p$ does not divide the dimension of $M$, where $p$ is the characteristic of $K$.*

*Proof.* $K$ is a direct summand of $M \otimes_K M^*$ if and only if $p$ does not divide the dimension of $M$, see [Ben, Proposition 3.1.9]. Thus, by assumption we have that $K$ is a direct summand of $M \otimes_K M^*$ and then the result follows by using $\mathrm{Ext}^i_{KG}(M,M) \cong \mathrm{Ext}^i_{KG}(K, M \otimes_K M^*)$, since $\mathrm{Ext}^i_{KG}(K, M \otimes_K M^*)$ is non-zero for all $i$ since it contains $\mathrm{Ext}^i_{KG}(K,K)$ as a direct summand and we assume that conjecture K is true. □

Recall that a finite group $G$ has a *normal $p$-complement* for a prime $p$ if $G$ is the semidirect product of a normal subgroup $H$ of order coprime to $p$ and a Sylow $p$-subgroup of $G$.

**Proposition 2.13.** *Let $G$ be a finite group with a normal $p$-complement and $K$ a field of characteristic $p$. If $S$ is a simple $KG$-module with $\mathrm{Ext}^1_{KG}(S,S) \neq 0$ then $\mathrm{Ext}^n_{KG}(S,S) \neq 0$ for all $n \geq 1$.*

*Proof.* Let $P_i$ denote the indecomposable projective $KG$-modules. Then $KG$ is isomorphic to the product of the algebras $M_{n_i}(\mathrm{End}_{KG}(P_i))$, the $n_i \times n_i$-matrix rings over the ring $\mathrm{End}_{KG}(P_i)$ for some integers $n_i$, see for example [HB, Theorem 14.9]. Thus the blocks of $KG$ are all Morita equivalent to the local algebras $\mathrm{End}_{KG}(P_i)$. Let

$$\cdots \to A^{m_i} \to \cdots A^{m_1} \to A^{m_0} \to S \to 0$$

be a minimal projective resolution of a simple non-projective module $S$. Then $\dim \mathrm{Ext}^i_A(S,S) = \dim \mathrm{Hom}_A(\Omega^i(S), S) = \dim \mathrm{Hom}_A(A^{m_i}, S) = m_i \neq 0$ for $i \geq 1$, by Lemma 3.1. □

## 3. Hybrid algebras and tame blocks of group algebras

We have already seen that it is not true for general symmetric algebras in Example 3.8 that $\mathrm{Ext}^1(S,S) \neq 0$ implies $\mathrm{Ext}^n(S,S) \neq 0$ for all $n \geq 1$. But we give a positive answer for certain simple modules in the class of hybrid algebras, as introduced in [ES2]. The class of hybrid algebras includes amongst others all Brauer graph algebras, and all weighted surface algebras (in particular almost all blocks of tame representation type). As a consequence, we get a positive answer for blocks of group algebras of tame representation type by considering extra cases at the end of this section. The results are as follows, the definitions are given below.

**Theorem 3.1.** *Assume $\Lambda$ is a hybrid algebra, where $\Lambda = KQ/I$. Let $i$ be a vertex of $Q$ so that for the simple module $S_i$ we have $\mathrm{Ext}^1(S_i,S_i) \neq 0$. Then $\mathrm{Ext}^n(S_i,S_i) \neq 0$ for all $n \geq 1$ unless $i$ is a hybrid vertex, and if $\beta$ is an arrow starting at $i$ which is not a loop then the module $\beta\Lambda$ does not have $\Omega$-period 3.*

We will show in the last section that the exception does not occur for blocks of group algebras. As well, we deal there with the potential blocks which are not hybrid algebras, and get:

**Corollary 3.2.** *Assume $\Lambda$ is a block of tame representation type with a simple module $S_i$ such that $\mathrm{Ext}^1(S_i,S_i) \neq 0$. Then $\mathrm{Ext}^n(S_i,S_i) \neq 0$ for all $n \geq 1$.*

Recall that for simple modules $S_i, S_j$ we have $\mathrm{Ext}^n_\Lambda(S_j,S_i) \cong \mathrm{Hom}_\Lambda(\Omega^n(S_j), S_i)$. This is non-zero if and only if the indecomposable projective module $P_i$ is a direct summand of the $n$-th term of a minimal projective resolution of $S_j$. Clearly, we may assume $\Lambda$ is not local.

Hybrid algebras are given by a presentation $\Lambda = KQ/I$, see below. We will first work with such algebras where $Q$ is the Gabriel quiver of $\Lambda$. After that, we will describe for a general hybrid algebra how to modify the proof. In the last section we give an example.

3.1. **Notation and Strategy.** For a module $M$ we write $\text{top}(M) := M/\text{rad}(M)$, and $P(M)$ is the projective cover of $M$. Consider an exact sequence

$$0 \to M_1 \to M \to M_2 \to 0$$

We call it *top good* if $\text{top}(M) = \text{top}(M_1) \oplus \text{top}(M_2)$. If so then the direct sum of projective covers of $M_1$ and of $M_2$ is a projective cover of $M$, and it follows that there is an exact sequence

$$0 \to \Omega(M_1) \to \Omega(M) \to \Omega(M_2) \to 0$$

Let $S_i$ be a simple module such that $\text{Ext}^1(S_i, S_i) \neq 0$, that is $\text{Hom}(\Omega(S_i), S_i) \neq 0$ so that the quiver has a loop at vertex $i$. For $\Lambda$ a hybrid algebra, the aim is to construct an exact sequence

$$0 \to W \to \Omega^2(S_i) \to S_i \to 0.$$

which is top good, and where applying $\Omega$ repeated times, it remains top good. This will then imply that $\text{Hom}(\Omega^n(S_i), S_i) \neq 0$ for all $n$.

In our case, $W$ will typically be a direct sum of modules of the form $\alpha\Lambda$ and moreover $\Omega(\alpha\Lambda) = \alpha_1\Lambda$ where $\alpha$ and $\alpha_1$ are arrows of $Q$. This is sufficient to have exact sequences which are top good. Namely if an arrow module $\alpha\Lambda$ is contained in some indecomposable non-projective module $M$ then $\alpha \notin \text{rad}(M)$, this follows from the fact that $\Lambda$ is symmetric hence selfinjective.

3.2. **Hybrid algebras.** We define the algebras in the generality as we need it: Results on Ext of simple modules do not depend on socle deformations, so we will not introduce these. As well, we start with algebras where $Q$ will be the Gabriel quiver of the algebra. Assume $(Q, f)$ is a biserial quiver. That is, $Q$ is a connected 2-regular quiver and $f$ is a permutation of the arrows of $Q$ such that for each arrow $\alpha$ we have $s(f(\alpha)) = t(\alpha)$, where $s(\beta)$ denotes the start of an arrow $\beta$ and $t(\beta)$ denotes the end (target) of the arrow $\beta$. We have an involution $(\overline{-})$ on the arrows, taking $\bar{\alpha}$ to be the arrow $\neq \alpha$ with the same starting vertex. Given $f$, this uniquely determines the permutation $g$ on arrows, defined by $g(\alpha) = \overline{f(\alpha)}$.

We fix a weight function $m_\bullet$ from $g$-cycles to $\mathbb{N}$, and a parameter function $c_\bullet$ from $g$-cycles to $K^*$. If $\alpha$ is an arrow, We write $m_\alpha$ and $c_\alpha$ for the values of its $g$-cycle, and $n_\alpha$ is the length of this cycle.

We define elements $A_\alpha$ and $B_\alpha$ in the path algebra $KQ$ as follows. We let $B_\alpha$ be the monomial along the $g$-cycle of $\alpha$ which starts with $\alpha$, of length $m_\alpha n_\alpha$, and $A_\alpha$ is the submonomial of $B_\alpha$ of length $m_\alpha n_\alpha - 1$ such that $B_\alpha = A_\alpha g^{-1}(\alpha)$. We write $|A_\alpha|$ for the length of $A_\alpha$.

3.2.1. **Notation.** The arrows in $f$-orbits of length 3 or 1 play a special role, we refer to these as *triangles.* In order to define a hybrid algebra with the above data, we fix a set $\mathcal{T}$ of triangles in $Q$, note that $\mathcal{T}$ is invariant under $f$.

*Assumption* 3.1. We assume that $m_\alpha n_\alpha \geq 3$ if $\bar{\alpha} \in \mathcal{T}$ and $m_\alpha n_\alpha \geq 2$ otherwise. That is, $|A_\alpha| \geq 2$ for $\alpha \in \mathcal{T}$ and $|A_\alpha| \geq 1$ otherwise. This means that $\alpha$ is always an arrow in the Gabriel quiver of the algebra.

Let $i$ be a vertex and $\alpha, \bar{\alpha}$ the arrows starting at $i$. We say that $i$ is *biserial* if $\alpha, \bar{\alpha}$ are both not in $\mathcal{T}$. We call the vertex $i$ a *quaternion* vertex if $\alpha, \bar{\alpha}$ are both in $\mathcal{T}$. Otherwise, we say that $i$ is *hybrid.*

**Definition 3.3.** The hybrid algebra $\Lambda = H_{\mathcal{T}} = H_{\mathcal{T}}(Q, f, m_\bullet, c_\bullet)$ associated to a set $\mathcal{T}$ of triangles is the algebra $\Lambda = KQ/I$ where $I$ is generated by the following elements:

(1) $\alpha f(\alpha) - c_{\bar{\alpha}} A_{\bar{\alpha}}$ for $\alpha \in \mathcal{T}$; and $\alpha f(\alpha)$ for $\alpha \notin \mathcal{T}$.
(2) $\alpha f(\alpha) g(f(\alpha))$,
(2') $\alpha g(\alpha) f(g(\alpha))$.
(3) $c_\alpha B_\alpha - c_{\bar{\alpha}} B_{\bar{\alpha}}$ for any arrow $\alpha$ of $Q$.

When $\mathcal{T} = \emptyset$, the algebra is a Brauer graph algebra. On the other extreme, when $\mathcal{T} = Q_1$ (and in particular $f^3 = \text{Id}$), the algebra is a weighted surface algebra.

3.2.2. *Some properties.* Suppose $\alpha, \bar{\alpha}$ are the arrows starting at vertex $i$.
(1) The arrows ending at vertex $i$ are $g^{-1}(\alpha)$ and $g^{-1}(\bar{\alpha})$. Moreover

$$g^{-1}(\alpha) = f^{-1}(\bar{\alpha}), \quad g^{-1}(\bar{\alpha}) = f^{-1}(\alpha)$$

(2) We have $\dim e_i\Lambda = m_\alpha n_\alpha + m_{\bar{\alpha}} n_{\bar{\alpha}}$. As well $m_{\bar{\alpha}} n_{\bar{\alpha}} = m_{f^{-1}\alpha} n_{f^{-1}\alpha}$ since $g^{-1}(\bar{\alpha}) = f^{-1}(\alpha)$, and $m, n$ are constant on $g$-cycles.
(3) If $\alpha \in \mathcal{T}$ then $\dim \alpha\Lambda = m_\alpha n_\alpha + 1$. Namely $\alpha\Lambda$ has basis consisting of the initial submonomials of $B_\alpha$ of length $\geq 1$ together with $A_{\bar{\alpha}}$. On the other hand if $\alpha \notin \mathcal{T}$ then $\dim \alpha\Lambda = m_\alpha n_\alpha$. In both cases, $\dim(\alpha g(\alpha)\Lambda) = m_\alpha n_\alpha - 1$.

3.3. **The simple module at a quaternion vertex.** Assume $i$ is a vertex such that both $\alpha, \bar{\alpha}$ are in $\mathcal{T}$. Then the simple module $S_i$ is periodic of period four, see [ES2]. There is an exact sequence which produces a minimal projective resultion, via slicing, of the form

$$0 \to S_i \to P_i \to P_i^- \to P_i^+ \to P_i \to S_i \to 0 \tag{1}$$

Here $P_i^+ = \oplus_{\alpha:i\to j} P_j$, that is, the summands are in bijection with the arrows starting at $i$. Similarly we have $P_i^- = \oplus_{\gamma:k\to i} P_k$ since it is the injective hull of $\Omega^{-1}(S_i)$. Then $P_i^+$ is a projective cover of $\Omega(S_i)$, and $P_i^-$ is a projective cover of $\Omega^2(S_i)$. Of the exceptions, the only one whose quiver has a loop is the singular disc algebra, which will be dealt with below in Lemma 4.11.

**Corollary 3.4.** *Assume there is a loop in the quiver at vertex $i$. Then* $\mathrm{Hom}(\Omega^n(S_i), S_i) \neq 0$ *for all* $n \geq 1$.

*Proof.* With the above, $P_i^+$ is a projective cover of $\Omega(S_i)$, and $P_i^-$ is a projective cover of $\Omega^2(S_i)$. □

3.4. **The simple module at a biserial vertex.** Assume $i$ is a vertex such that $\alpha$ and $\bar{\alpha}$ are both not in $\mathcal{T}$. Write $\alpha_j = f^j(\alpha)$ and $\bar{\alpha}_j = f^j(\bar{\alpha})$.

We identify $\Omega(S_i) = \alpha\Lambda + \bar{\alpha}\Lambda$, and it has projective cover $P_{t\alpha} \oplus P_{t(\bar{\alpha})} \xrightarrow{d_1} \alpha\Lambda + \bar{\alpha}\Lambda$ where $d_1(x, y) = \alpha x + \bar{\alpha} y$.

We identify $\ker(d_1)$ with $\Omega^2(S_i)$. This contains $(\alpha_1, 0)\Lambda$ and also $(0, \bar{\alpha}_1)\Lambda$ since $\alpha\alpha_1 = 0$ and $\bar{\alpha}\bar{\alpha}_1 = 0$. These generate submodules of $\Omega^2(S_i)$ which are not contained in the radical. Furthermore, the relation (3) of Definition 4.3 shows that $\Omega^2(S_i)$ also contains

$$\psi := (c_\alpha A_{g\alpha}, -c_{\bar{\alpha}} A_{g\bar{\alpha}}).$$

This does not lie in the submodule $W := (\alpha_1, 0)\Lambda \oplus (0, \bar{\alpha}_1)\Lambda$ of $\Omega^2(S_i)$. In fact $\psi\Lambda$ intersects each of the two summands of $W$ in its socle, and $\psi\Lambda/W \cap \psi\Lambda$ is simple and isomorphic to $S_i$. By dimensions, $W + \psi\Lambda = \Omega^2(S_i)$.

This shows that we have an exact sequence

$$0 \to W \to \Omega^2(S_i) \to S_i \to 0 \tag{2}$$

and moreover it is top good.

**Lemma 3.5.** *For each $r \geq 0$ we have an exact sequence*

$$0 \to \Omega^r(W) \to \Omega^{r+2}(S_i) \to \Omega^r(S_i) \to 0. \tag{3}$$

*Proof.* The case $r = 0$ is discussed above. Since $W$ is a direct sum of modules generated by arrows, the sequence for $r = 0$ is top good (see subsection 4.1). Hence we get the sequence for $r = 1$. Now, $\Omega(W)$ is isomorphic to the direct sum $\alpha_2\Lambda \oplus \bar{\alpha}_2\Lambda$ and hence the sequence for $r = 1$ also is top good and we get the sequence for $r = 2$. The lemma follows by induction on $r$. □

**Corollary 3.6.** *Assume there is a loop at vertex $i$. Then* $\mathrm{Hom}(\Omega^n(S_i), S_i) \neq 0$ *for all* $n \geq 1$.

*Proof.* Say, $\alpha$ is a loop, we have an exact sequence

$$0 \to \bar{\alpha}\Lambda \to \Omega(S_i) \to \alpha\Lambda/\langle B_\alpha \rangle \to 0$$

which implies $0 \neq \mathrm{Hom}(\alpha\Lambda/\langle B_\alpha\rangle, S_i) \subseteq \mathrm{Hom}(\Omega(S_i), S_i)$. The sequence (3) for $r = 0$ shows $\mathrm{Hom}(\Omega^2(S_i), S_i) \neq 0$, and then by induction the Corollary follows from the Lemma. □

3.5. **The simple module at hybrid vertex.** Let $\alpha$ be an arrow in an $f$-cycle of length 3 or 1, and set $i_j = s(\alpha_j)$. Assume $\alpha \in \mathcal{T}$ but $\bar{\alpha}$ is not in $\mathcal{T}$. Write $\alpha = \alpha_0$ and $\alpha_r = f^r(\alpha)$ (possibly $\alpha_r = \alpha$ for all $r$). Let also $\beta = \beta_0 = \bar{\alpha}$ and $\beta_j = f^j(\beta)$. Therefore $\alpha\alpha_1 = c_\beta A_\beta$ and $\beta\beta_1 = 0$. As well, $g(\alpha_2) = \bar{\alpha} = \beta$.

We have $\Omega(S_i) = \alpha\Lambda + \beta\Lambda$, and $\beta\Lambda$ is periodic but $\alpha\Lambda$ is not. It is natural to use the exact sequence

$$0 \to \beta\Lambda \to \Omega(S_i) \to \alpha\Lambda/(\alpha\Lambda \cap \beta\Lambda) \to 0.$$

Therefore we need a good description of the first few syzygies of the cokernel. Note that for the moment we do not assume that $\alpha$ or $\bar{\alpha}$ is a loop.

**Lemma 3.7.** *Consider the first few syzygies of the the module $\alpha\Lambda/(\alpha\Lambda \cap \beta\Lambda)$. We have*
*(a)* $\Omega(\alpha\Lambda/(\alpha\Lambda \cap \beta\Lambda)) \cong \alpha_1\Lambda$.
*(b)* $\Omega(\alpha_1\Lambda) \cong \alpha_2\beta\Lambda$.
*(c) There is an exact sequence* $0 \to \beta_1\Lambda \to \Omega(\alpha_2\beta\Lambda) \to S_{i_0} \to 0$

We will prove this below. First we will show how this is applied. The Lemma gives exact sequences for the first few syzygies, where the kernels are generated by arrows.

$$\begin{aligned}
0 \longrightarrow \beta\Lambda \longrightarrow \Omega(S_{i_0}) \longrightarrow \alpha\Lambda/(\alpha\Lambda \cap \beta\Lambda) \longrightarrow 0\\
0 \longrightarrow \beta_1\Lambda \longrightarrow \Omega^2(S_{i_0}) \longrightarrow \alpha_1\Lambda \longrightarrow 0\\
0 \longrightarrow \beta_2\Lambda \longrightarrow \Omega^3(S_{i_0}) \longrightarrow \alpha_2\beta\Lambda \longrightarrow 0\\
0 \longrightarrow \beta_3\Lambda \longrightarrow \Omega^4(S_{i_0}) \longrightarrow U = \Omega((\alpha_2\beta)\Lambda) \longrightarrow 0\\
0 \longrightarrow \beta_1\Lambda \longrightarrow U \longrightarrow S_{i_0} \to 0
\end{aligned}$$

The sequence for $\Omega^4(S_{i_0})$ can be rearranged to

$$0 \longrightarrow W' \longrightarrow \Omega^4(S_{i_0}) \longrightarrow S_{i_0} \longrightarrow 0$$

where $W'$ is an extension $0 \to \beta_3\Lambda \to W' \to \beta_1\Lambda \to 0$. All sequences are top good, and the syzygies of the kernels are also top good. With this, we get the following.

**Corollary 3.8.** *(a) If $\beta$ is a loop then* $\mathrm{Hom}(\Omega^n(S_{i_0}), S_{i_0}) \neq 0$ *for all $n \geq 1$.*
*(b) Suppose $\alpha$ is a loop. Then* $\mathrm{Hom}(\Omega^n(S_{i_0}), S_{i_0}) \neq 0$ *for all $n \not\equiv 3$ modulo 4.*
*(c) If in case (b) the arrow $\beta$ lies in a 3-cycle of $f$ then* $\mathrm{Hom}(\Omega^n(S_{i_0}), S_{i_0}) \neq 0$ *for any $n \geq 1$.*

*Proof.* (a) The algebra is not local, so if $\beta$ is a loop then $\alpha$ is not a loop and belongs to a 3-cycle of $f$. Then we see that $f(\beta) = \beta$, that is $\beta_i = \beta$ for all $i$. Since the sequences are top good and $\mathrm{Hom}(\beta_i\Lambda, S_{i_0}) \neq 0$ part (a) follows.

(b) Here necessarily $f(\alpha) = \alpha$ and $\alpha_i = \alpha$ for all $i$. Therefore $\mathrm{Hom}(V, S_{i_0}) \neq 0$ for $V$ one of $\alpha\Lambda/\alpha\Lambda\cap\beta\Lambda$, or $\alpha_1\Lambda$, or $U$, *but not for* $V = \alpha_2\beta\Lambda$. This implies part (b).

(c) With this assumption we have $\beta_2$ ends at vertex $i_0$, therefore $\mathrm{Hom}(\beta_2\Lambda, S_{i_0}) \neq 0$. The third exact sequence above is top good and it follows that $\mathrm{Hom}(\Omega^3(S_{i_0}), S_{i_0}) \neq 0$. Apply $\Omega^3$ to the sequence for $\Omega^4(S_{i_0})$, this gives

$$0 \to \Omega^3(W') \to \Omega^7(S_{i_0}) \to \Omega^3(S_{i_0}) \to 0$$

Since $\mathrm{Hom}(\Omega^3(S_{i_0}), S_{i_0}) \neq 0$ we deduce that $\mathrm{Hom}(\Omega^7(S_{i_0}), S_{i_0}) \neq 0$. Part (c) follows by induction. □

**Proof of Lemma 4.7** (a) The minimal relation $\alpha\alpha_1 = c_\beta A_\beta$ shows that $\alpha\Lambda \cap \beta\Lambda = (\alpha\alpha_1)\Lambda$, and moreover it is 2-dimensional (we know $\alpha\alpha_1\alpha_2$ spans the socle and $\alpha\alpha_1 g(\alpha_1) = 0$).

Hence if $\varphi : e_{i_1}\Lambda \to \alpha\Lambda/(\alpha\alpha_1\Lambda)$ is the composite of the map $x \to \alpha x$ with the canonical surjection, then we have $\varphi(\alpha_1\Lambda) = 0$. To show that the kernel of $\varphi$ is equal to $\alpha_1\Lambda$, we need that $\dim \alpha_1\Lambda + \dim(\alpha\Lambda/(\alpha\Lambda \cap \beta\Lambda)) = \dim e_{i_1}\Lambda$ First, $\dim \alpha_1\Lambda = m_{\alpha_1}n_{\alpha_1} + 1$, and $\dim(\alpha\Lambda/(\alpha\alpha_1\Lambda)) = m_\alpha n_\alpha - 1$. The sum of these is

$$m_{\alpha_1}n_{\alpha_1} + m_\alpha n_\alpha = \ m_{\alpha_1}n_{\alpha_1} + m_{\bar{\alpha}_1}n_{\bar{\alpha}_1} = \dim e_{i_1}\Lambda$$

(b) By the relations, $\alpha_1\alpha_2 g(\alpha_2) = 0$ which shows that $\alpha_2 g(\alpha_2)\Lambda \subseteq \Omega(\alpha_1\Lambda)$. We have

$$\dim(\alpha_2 g(\alpha_2)\Lambda) = m_{\alpha_2}n_{\alpha_2} - 1, \quad \dim(\alpha_1\Lambda) = m_{\alpha_1}n_{\alpha_1} + 1$$

and the sum is equal to

$$m_{\alpha_1}n_{\alpha_1} + m_{\alpha_2}n_{\alpha_2} = m_{\bar{\alpha}_2}n_{\bar{\alpha}_2} + m_{\alpha_2}n_{\alpha_2} = \dim(e_{i_2}\Lambda)$$

noting that $\bar{\alpha}_2 = g(\alpha_1)$. Hence $\Omega(\alpha_1\Lambda) = \alpha_2\beta\Lambda$, as stated.

(c) Let $U := \Omega(\alpha_2\beta\Lambda)$. Let $x = t(\beta)$, then

$$\dim U = \dim e_x\Lambda - m_{\alpha_2}n_{\alpha_2} + 1 = m_{\beta_1}n_{\beta_1} + m_{\bar{\beta}_1}n_{\bar{\beta}_1} - m_{\alpha_2}n_{\alpha_2} + 1.$$

Since $g(\alpha_2) = \beta$ and $g(\beta) = \bar{\beta}_1$, the corresponding terms cancel and we get $\dim U = m_{\beta_1}n_{\beta_1} + 1$. Now, $\alpha_2\beta\beta_1 = 0$ which shows that $\beta_1\Lambda$ is contained in $U$; it has dimension $m_{\beta_1}n_{\beta_1}$. Furthermore, we have

$$\alpha_2\beta A_{\bar{\beta}_1} = \alpha_2 B_\beta = 0$$

Now, $A_{\bar{\beta}_1} \notin \beta_1\Lambda$ and it follows that

$$U = \beta_1\Lambda + A_{\bar{\beta}_1}\Lambda.$$

Moreover $U/\beta_1\Lambda$ is simple, and it is isomorphic to $S_{i_0}$: We know $\beta A_{\bar{\beta}_1}$ is a socle element and the algebra is symmetric, therefore also $A_{\bar{\beta}_1}\beta$ is a socle element which shows that $A_{\bar{\beta}_1} = A_{\bar{\beta}_1}e_{i_0}$. This gives the stated exact sequence. □

3.6. **Hybrid algebras with virtual arrows.** The definition of hybrid algebras where $Q$ may not be the Gabriel quiver, is slightly more general:

First, replace Assumption 4.1 by a weaker assumption. Namely assume only that $|A_\alpha| \geq 1$ for $\bar{\alpha} \in \mathcal{T}$, and that $|A_\alpha| \geq 0$ for $\alpha, \bar{\alpha} \notin \mathcal{T}$.

**Definition 3.9.** We call an arrow $\alpha$ *virtual* if one of the following holds.
(i) $|A_\alpha| = 1$ and $\bar{\alpha} \in \mathcal{T}$, or
(ii) $|A_\alpha| = 0$ and $\alpha, \bar{\alpha} \notin \mathcal{T}$.

This means that $\alpha$ is virtual if and only if it is not in the Gabriel quiver.

Next, replace the zero relations (2) and (2') in Definition 4.3 by the following.

(2)* $\alpha f(\alpha)g(f(\alpha)) = 0$ in $\Lambda$ unless $\alpha, \bar{\alpha} \in \mathcal{T}$ and $|A_{\bar{\alpha}}| = 1$, or else $\alpha, \bar{\alpha}, g(\alpha) \in \mathcal{T}$ and $|A_{\bar{\alpha}}| = 2$ and $|A_{f(\alpha)}| = 1$.

(2')* $\alpha g(\alpha)f(g(\alpha)) = 0$ in $\Lambda$ unless $\alpha, g(\alpha) \in \mathcal{T}$ and $|A_{f(\alpha)}| = 1$, or else $g(\alpha), \alpha, \bar{\alpha} \in \mathcal{T}$ and $|A_{f(\alpha)}| = 2$ and $|A_{\bar{\alpha}}| = 1$.

Here we do not need to use these directly. We describe now the modifications needed. In the last section we give an example.

We assume always that there is a loop at vertex $i$ which is part of the Gabriel quiver. Moreover, we assume that $Q$ is connected with at least two vertices. This implies that the arrow $\bar{\alpha}$ cannot be virtual. After the changes to be made, the results remain the same. We give an example later.

Suppose we have $r \geq 1$ and an exact sequence

$$0 \to \Omega^r(W) \to \Omega^{r+2}(S_{i_0}) \to \Omega^r(S_{i_0}) \to 0$$

which is not top good. This happens if and only if the module $\Omega^r(W)$ has a direct summand $\gamma\Lambda$ with $\gamma$ a virtual arrow. Then we get from this an exact sequence

$$0 \to \Omega^{r+1}(W) \to P' \oplus \Omega^{r+3}(S_{i_0}) \to \Omega^{r+1}(S_{i_0}) \to 0$$

and $P' = P(\gamma\Lambda)$. This cannot be $P(S_{i_0})$ since there is no virtual arrow ending at vertex $i_0$. So this does not have effect on $\mathrm{Hom}(-, S_{i_0})$. We continue as before.

Note that virtual arrows along the $f$-cycle of arrows not in $\mathcal{T}$ are isolated, that is if $\gamma$ is virtual then $f(\gamma)$ is not virtual.

3.7. **Tame blocks.** All tame blocks occur amongst the algebras of dihedral, semidihedral or quaternion type. Almost all of these algebras are hybrid algebras, with a few exceptions. We explain first why Conjecture S holds for hybrid algebras, and then we will deal with the exceptions.

Let $\Lambda$ be an algebra of dihedral, semidihedral or quaternion type, and let $i$ be a vertex of the quiver with a loop.

(1) Assume $\Lambda$ is a hybrid algebra. If $i$ is biserial, or quaternion, then the claim follows from Corollary 4.4 and 4.6. Suppose $i$ is a hybrid vertex, then $\Lambda$ is of semidihedral type. In Corollary 4.8 we had the exception when $n \equiv 3$ modulo 4 which occurs when $\alpha$ is a loop in $\mathcal{T}$. In that case, for the algebras of semidihedral type, the arrow $\bar{\alpha}$ always belongs to an $f$-cycle of length three. Therefore we get also that $\mathrm{Hom}(\Omega^{4k+3}(S_i), S_i) \neq 0$ for all $k \geq 0$. There is one exception; namely, the algebra $SD(2\mathcal{B})_3$ with $s = 2$ but this will be dealt with below in Lemma 4.11.

(2) We consider now algebras of quaternion type which are not hybrid algebras. They are the algebras $Q(3\mathcal{C})$ and $Q(2\mathcal{B})_2$ in [E]. In both cases, each simple module has $\Omega$-period four (this follows since the algebras are periodic of period four, see [ES4, Theorem 5.1], based on [H]). Therefore, there is an exact sequence of the form (1), and the arguments for Corollary 4.4 can be used here as well, and the claim follows.

(3) This leaves the case when $\Lambda$ is of semidihedral type, listed as $SD(3\mathcal{C})$, or $SD(2\mathcal{B})_3$ in [E], they are not hybrid algebras. We show that in each of these cases if $S_0$ is the simple module at a loop then $\mathrm{Hom}(\Omega^n(S_0), S_0) \neq 0$ for all $n \geq 0$, this will be done in separate lemmas. Assuming these, the proposition will follow.

We start with the case when $\Lambda$ is of type $SD(3\mathcal{C})$. Explicitly, $SD(3\mathcal{C})_1$ is defined as $SD(3\mathcal{C})_1 = kQ/I_1$, where $Q$ is the quiver

$1 \underset{\gamma}{\overset{\beta}{\rightleftarrows}} 0 \underset{\eta}{\overset{\delta}{\rightleftarrows}} 2$ (with a loop $\rho$ at $0$) and $I_1 = (\beta\delta, \beta\rho, \rho\gamma, \eta\gamma, \eta\rho, \rho\delta, \rho^s - \gamma\beta, \gamma\beta - \delta\eta, \beta\gamma\beta, \eta\delta\eta)$ with $s \geq 3$.

Moreover, $SD(3\mathcal{C})_2$ is defined as $SD(3\mathcal{C})_2 = kQ/I_2$ with $I_2 = (\beta\rho, \rho\delta, \eta\rho, \rho\gamma, \gamma\beta - \delta\eta, (\gamma\beta)^k - \rho^s, (\beta\gamma)^{k-1}\beta\gamma, (\eta\delta)^{k-1}\eta\gamma)$ where $s \geq 2, k \geq 2$. We remark that the relation $(\gamma\beta)^k - \rho^s$ is needed instead of the relation $(\beta\gamma)^k - \rho^s$, which is a typo in [E] on page 301.

**Lemma 3.10.** *Assume $\Lambda$ is of type $SD(3\mathcal{C})_1$, or $SD(3\mathcal{C})_2$. Then there is an exact sequence*

$$0 \to W \to \Omega^2(S_0) \to S_0 \to 0$$

*which is top good, with $W$ periodic of period 3, such that taking syzygies of the sequence gives rise to top good sequences.*

*Proof.* For this algebra, we use syzygies of $\rho\Lambda$.

Assume first we have type $SD(3\mathcal{C})_1$. Then $U := \Omega(\rho\Lambda)$ has Loewy length two, its socle is $S_0$ and its top is $S_1 \oplus S_2$. From this, we see $\Omega(U)$ has Loewy length two, it has simple top $S_0$ and socle isomorphic to $S_1 \oplus S_2$, and $\Omega^2(U) \cong \rho\Lambda$.

One checks that there is an exact sequence, which is top good,

$$0 \to W \to \Omega^2(S_0) \to S_0 \to 0$$

where $W$ is an extension of the form

$$0 \to U \to W \to \Omega(U) \to 0$$

and $W$ is periodic, also of period 3.

Now assume we have $SD(3\mathcal{C})_2$. Then $\rho\Lambda$ also has $\Omega$-period three. The syzygies are slightly more complicated. First, $U := \Omega(\rho\Lambda)$ can be constructed as follows. Take the direct sum $P_1 \oplus P_2$ and note that $\mathrm{rad}(P_1)/S_1$ is isomorphic to $\mathrm{rad}(P_2)/S_2$. Therefore if we factor out the submodule $(\beta, \eta)\Lambda$ then the quotient has radical isomorphic to $\mathrm{rad}(P_1)/S_1$ and has top $S_1 \oplus S_2$.

Then $\Omega(U)$ has socle isomorphic to $S_1 \oplus S_2$ and $\Omega(U)/S_1 \oplus S_2$ is isomorphic to $\mathrm{rad}(P_1)/S_1$. Furthermore, $\Omega^2(U) \cong \rho\Lambda$. Then the result follows as in the previous case. □

The only other possibility is that $\Lambda$ is the algebra $SD(2\mathcal{B})_3$. This is the same as the singular disc algebra introduced in [ES5], which was mentioned as an exception in section 4.3. This algebra has a better presentation, which can be found in [H, Proposition 4.2] (where it is called $SD(2\mathcal{B})_4^s(c)$), which we use here. Let $Q$ be the quiver

$$\alpha \circlearrowright 0 \underset{\gamma}{\overset{\beta}{\rightleftarrows}} 1 \circlearrowleft \eta$$

Then $\Lambda = KQ/I$ where $I$ is generated by

$$\beta\gamma = \alpha^2,\ \alpha\beta = \beta\eta,\ \eta\gamma = \gamma\alpha,\ \gamma\beta = \eta^2(1 + c\eta^{s-1}),\ (s \geq 2)$$

for some $c \in K$, and the zero relations $\alpha^s\beta = 0 = \alpha^{s+2} = \eta^{s+2} = \eta^s\gamma = \beta\eta^s$. Note that the relations are essentially symmetric if the vertices are interchanged. Namely, we can replace $\beta$ by $\beta' := \beta u^{-1}$ where $u$ is the unit $u = 1 + c\eta^{s-1}$. Observing also that $\beta\eta^{s-1}\gamma$ is a scalar multiple of $\alpha^{s+1}$, the relations become

$$\beta'\gamma = \alpha^2(1 + d\alpha^{s-1}),\ \alpha\beta' = \beta'\eta,\ \eta\gamma = \gamma\alpha,\ \gamma\beta' = \eta^2\ (d \in K)$$

with the same zero relations. Therefore, it suffices to analyse the syzygies for $S_0$.

The exceptional property of the algebra is that $\mathrm{rad}P_0/S_0$ is isomorphic to $\mathrm{rad}P_1/S_1$. We denote this module by $H$. The algebra $\Lambda$ has Cartan matrix $\begin{pmatrix} s+2 & s \\ s & s+2 \end{pmatrix}$. Note that $H$ has dimension vector $(s, s)$. The socle of the algebra is spanned by $\alpha^{s+1}$ and $\beta^{s+1}$.

When $s = 2$, this is a hybrid algebra. Both vertices are quaternion, however the simples are not $\Omega$-periodic.

**Lemma 3.11.** *Assume $\Lambda$ is of type $SD(2\mathcal{B})_3$. Then there is a top good exact sequence*

$$0 \to W \to \Omega^2(S_0) \to S_1 \to 0$$

*with $W \cong \Omega^2(W)$. Moreover,* $\mathrm{Hom}(\Omega^r(W), S_i) \neq 0$ *for $i = 0, 1$ and $r \geq 0$.*

*Proof.* We start with the projective cover $\pi : P_0 \oplus P_1 \to \alpha\Lambda + \beta\Lambda \cong \Omega(S_0)$, which is given by $\pi(x, y) = \alpha x + \beta y$, and we take $\Omega^2(S_0) = \ker(\pi)$. The relations show that $\Omega^2(S_0)$ contains the submodule

$$W := (\alpha, -\gamma)\Lambda + (\beta, -\eta)\Lambda$$

We will show that $W$ is a maximal submodule of $\Omega^2(S_0)$, with quotient the simple module $S_1$. Let $\phi = (\alpha, -\gamma)$ and $\psi := (\beta, -\eta)$, the two independent generators of $W$ and let $J$ be the radical of $\Lambda$.

(i) We claim that $\phi\Lambda \cap \psi\Lambda = \phi J = \psi J$.

Clearly, $\phi$ and $\psi$ are not in $\phi\Lambda \cap \psi\Lambda$. So it suffices to show that $\phi J = \psi J$. We have $\phi J = \langle \phi\alpha, \phi\beta \rangle$. Directly, $\phi\alpha = \psi\gamma$. Furthermore,

$$\phi\beta = \psi\eta + (0, -c\eta^{s+1}) = \psi\eta(1 - c\eta^s),$$

which is equal to $\psi\eta$ times a unit. It follows that $\phi J = \langle \psi\gamma, \psi\eta \rangle = \psi J$.

We have $\phi\Lambda \cong \Omega^{-1}(S_0)$ and hence $WJ = \phi J$ is isomorphic to $H$, which shows that $\underline{\dim}\, W = (s+1, s+1)$. On the other hand, we know $\underline{\dim}\, \Omega^2(S_0) = \underline{\dim}\, P_1 + \underline{\dim}\, S_0 = (s+1, s+2)$, and therefore, the quotient $\Omega^2(S_0)/W$ is simple and isomorphic to $S_1$. Explicitly, we can write down an element in the kernel of $\pi$ which is not in $W$, such as $(\beta\eta^{s-1}, 0) \in \Omega^2(S_0)$. We get the exact sequence as stated. This is top good since $W$ has maximal Loewy length. We also see from the generators that $\mathrm{Hom}(W, S_i) \neq 0$. We claim that $\Omega^2(W) \cong W$.

(i) Consider a projective cover $d : P_0 \oplus P_1 \to W$ given by $d(x, y) := \phi x + \psi y$. We compute (directly) $d(\alpha, -\gamma) = 0$ and $d(\beta, -\eta(1 + c\eta^{s-1})) = 0$. Hence $\Omega(W)$ contains the submodule $W' := \phi_1\Lambda + \psi_1\Lambda$ where $\phi_1 = (\alpha, -\gamma)$ and $\psi_1 = (\beta, -\eta(1 + c\eta^{s-1}))$. As before we compute that the dimension vector of $W'$ is $(s + 1, s + 1)$ which is the same as $\underline{\dim}\, \Omega(W)$. Hence $W' = \Omega(W)$, and also $\mathrm{Hom}(W', S_i) \neq 0$.

(ii) Consider a projective cover $d_1 : P_0 \oplus P_1 \to \Omega(W)$, where $d_1(x, y) = \phi_1 x + \psi_1 y$. We compute that the kernel of $d_1$ contains $\phi$ and $\psi$ and hence has submodule $W$. By dimensions, it is equal to $W$. □

The sequence in the previous lemma is top good and since $\Omega(W)$ also has maximal Loewy length, the sequence $0 \to \Omega(W) \to \Omega^3(S_0) \to \Omega(S_1) \to 0$ is also top good. Hence, applying $\mathrm{Hom}(-, S_0)$ takes any $0 \to \Omega^r(W) \to \Omega^{r+2}(S_0) \to \Omega^r(S_1) \to 0$ to an exact sequence, and therefore, $\mathrm{Hom}(\Omega^{r+2}(S_0), S_0) \neq 0$ since $\mathrm{Hom}(\Omega^r(W), S_0) \neq 0$ for every $r \geq 0$.

**Corollary 3.12.** *Assume $\Lambda$ is an algebra $SD(3\mathcal{C})_i$ for $i = 1$ or $2$, or $SD(2\mathcal{B})_3$. Then for all $n \geq 1$ we have* $\mathrm{Hom}(\Omega^n(S_0), S_0) \neq 0$. *For $SD(2\mathcal{B})_3$, we have as well that* $\mathrm{Hom}(\Omega^n(S_1), S_1) \neq 0$ *for all $n \geq 1$, by the symmetry of the relations.*

**Corollary 3.13.** *Assume $\Lambda$ is a tame block. If $S$ is simple and* $\mathrm{Ext}^1_\Lambda(S, S) \neq 0$ *then* $\mathrm{Ext}^n_\Lambda(S, S) \neq 0$ *for all $n \geq 1$.*

*Proof* This follows from Corollary 4.4, 4.6, 4.8, 4.12 and (1), (2) and (3) of Section 4.7.

**Example 3.14.** We describe the algebras of type $SD(2\mathcal{A})_2$. This has quiver

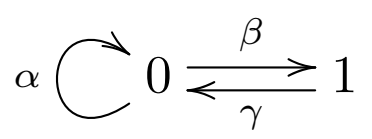

The relations are

$$\gamma\beta = 0,\ \alpha^2 = (\beta\gamma\alpha)^{k-1}\beta\gamma + c(\beta\gamma\alpha)^k,\ \ (\alpha\beta\gamma)^k = (\beta\gamma\alpha)^k.$$

Here $c \in K$ and $k \geq 2$. We do not need to consider socle deformation, that is we can take $c = 0$. This is a hybrid algebra where the quiver $Q$ is the above quiver extended by a virtual loop, call it $\varepsilon$, at vertex 1. We have

$$f = (\alpha)(\beta\ \varepsilon\ \gamma),\ \ g = (\alpha\ \beta\ \gamma)(\varepsilon)$$

In particular the relations involving $\varepsilon$ are

$$\varepsilon = (\gamma\alpha\beta)^k,\ \beta\varepsilon = 0,\ \varepsilon\gamma = 0.$$

With the notation as in Lemma 4.7 we have $\alpha_1 = \alpha$ and $\beta_1 = \varepsilon$ and $\beta_2 = \gamma$.

Let $k = 2$, then $X := \alpha\Lambda/\alpha\Lambda \cap \beta\Lambda$ is uniserial of the form $\mathcal{U}(0,1,0,0,1)$ (describing a uniserial module by its composition factors, from top to socle). Next, $\Omega(X) = \alpha\Lambda$ and $\Omega^2(X) = (\alpha\beta)\Lambda$. Furthermore, $\Omega^3(X)$ has length 2 with socle $S_1$ and top $S_0$.

Here $\Omega(\beta\Lambda) = S_1$ and we should interpret $\beta_1$ as the virtual loop $\varepsilon$. The exact sequence

$$0 \to S_1 \to \Omega^2(S_0) \to \alpha\Lambda \to 0 \tag{1}$$

is not top good. We get

$$0 \to \Omega(S_1) = \gamma\Lambda \to P_1 \oplus \Omega^3(S_0) \to (\alpha\beta)\Lambda \to 0 \tag{2}$$

Constructing it directly, we see $\Omega^3(S_0) = (\alpha\beta\gamma)\Lambda = \mathrm{rad}(\alpha\beta)\Lambda$ which then is consistent with the short exact sequence (2).

The sequence (2) is top good and we get

$$0 \to \Omega(\gamma\Lambda) = \beta\Lambda \to \Omega^4(S_0) \to \Omega^3(X) \to 0$$

We may rearrange this as in Lemma 4.7. In this case the module $W$ is an extension of the simple module $S_1$ by $\beta\Lambda$. We can verify directly that $W$ is indecomposable and periodic of period 3.

## 4. Appendix by Bernhard Böhmler and René Marczinzik : Magma code

In this appendix we present the Magma code that was used to verify the the following conjecture, that we call conjecture $S_n$ in the following:

**Conjecture $S_n$ 4.1.** *Let $n \geq 2$ be a fixed integer. Let $S$ be a simple $KG$ module with $\mathrm{Ext}^1_{KG}(S,S) \neq 0$. Then $\mathrm{Ext}^i_{KG}(S,S) \neq 0$ for all $i = 1, 2, ..., n$.*

The main result of this appendix can be stated as follows:

**Theorem 4.2.**
*(1) Conjecture $S_5$ is true for every group algebra $KG$ with $|G| \leq 500$, where $K$ is a splitting field of $G$ of characteristic $p$, where $p$ divides the group order of $G$.*
*(2) Conjecture $S_3$ is true for every group algebra $KG$ with $|G| \leq 1000$, where $K$ is a splitting field of $G$ of characteristic $p$, where $p$ divides the group order of $G$.*

The theorem was verified using the computer algebra system Magma and the following code that can be copied directly into Magma:

```
// The input of the function NthSyzygy is given by a KG-module M and a non-negative integer n.
// The group G and the field K are assumed finite.
// The function computes and returns the n-th syzygy of M.

NthSyzygy := function(M,n)
    if IsZero(n) then
        return M;
    end if;
    for i in [1..n] do
        P,f := ProjectiveCover(M);
        OMEGA := Kernel(f);
        M:=OMEGA;
    end for;
    return M;
end function;
```

```
// The input of the function DimensionExt_i_M_M is given by a KG-module M and a non-negative integer i.
// The group G and the field K are assumed finite.
// The function computes and returns the K-dimension of the i-th Ext group between M and M.

DimensionExt_i_M_M := function(M,i) // The variable i is assumed larger than 1 here.
    return Dimension(Ext(NthSyzygy(M,i-1),M));
end function;



// The input of the function ExtTestAuxiliary is given by G, p, and m, where G denotes
// a finite group, p denotes a prime number dividing the order of G, and m denotes a
// positive integer larger than or equal to 2.
// This function tests if for a given group G with minimal splitting field K over the
// field with p elements the following is satisfied: for every simple KG-module S where
// Ext_{KG}^1(S,S) is non-zero, also Ext_{KG}^i(S,S) is non-zero where 0<i<m+1.
// The function ExtTestAuxiliary returns either a counter-example or the boolean true.

ExtTestAuxiliary := function(G,p,m) // The integer m is assumed greater than 2 here.
    SIMS := AbsolutelyIrreducibleModules(G,GF(p));
    temp_field_sizes:=[];
    for S in SIMS do
        Append(~temp_field_sizes,#BaseRing(S));
    end for;
    ExponentsFieldOrders := [];
    for g in temp_field_sizes do
Append(~ExponentsFieldOrders, Floor(Log(p, g)));
    end for;
    LCM := Lcm(ExponentsFieldOrders);
    K := GF(p^LCM);
    Y := [ ExtendField(M, K) : M in SIMS ];
    s := #Y;
    tempext1 := [];
    j := 1; // This variable runs from 1 to m-1.
    for z in [1..s] do
        printf "During the computation of tempext_%o the variabe z is at the moment equal to %o. \n", j, z;
        M:=Y[z];
        if Dimension(Ext(M,M)) gt 0 then
            Append(~tempext1,M);
        end if;
    end for;

    tempext_old := tempext1;
    tempext_new := [];
    t := #tempext_old;

    while j lt m do

        for z in [1..t] do
            printf "During the computation of tempext_%o the variabe z is at the moment equal to %o. \n", j+1, z;
            N:=tempext_old[z];
            if DimensionExt_i_M_M(N,j+1) gt 0 then
                Append(~tempext_new, N);
            end if;
        end for;

        u := #tempext_new;

        if (IsZero(t-u) eq false) then
            print("We have found a counter-example: "); print([*G,p,j+1*]);
            return([*G,p,j+1*]);
        end if;

        j := j+1;
        tempext_old := tempext_new;
        tempext_new := [];
        t := #tempext_old;
    end while;

    return true;
end function;



// For a given group G and a prime number p dividing the order of G,
// the auxiliary funtion HasNormalpComplementSpecialCase does
// the following:
// if the order of G is of the form p^s*q^t for another prime number q,
// then it tests whether G has a normal p-complement (which is equivalent to
// the q-Sylow subgroup being normal). If the order of G is not of the above form
// then the boolean false is returned.

HasNormalpComplementSpecialCase := function(G,p)

    flag := false;

    primedivs := PrimeDivisors(#G);
    if #primedivs eq 2 then
        L := [x:x in primedivs | x ne p]; q := L[1];
        if IsNormal(G, SylowSubgroup(G,q)) then
            flag := true;
        end if;
```

```
    end if;

    return flag;
end function;



// The following is the main program. The input is given by three positive integers N, a, and b.
// The program tests for each group G whose order lies in the interval [a,b] if Ext_{KG}^i(S,S) is
// non-zero for all 0<i<N+1 provided that Ext_{KG}^1(S,S) is non-zero.
// Here, K denotes a minimal splitting field of G over a finite field with p elements
// and p runs through all prime divisors of the order of G.

// When the conjecture holds for all groups in the given interval and for all splitting fields,
// then the boolean true is returned; otherwise the boolean false is returned.

ExtTestSimples := function(N,a,b)

    Complete_List_boolean_for_entered_groups_and_primes := [];
    for n in [a..b] do
        P := SmallGroupProcess(n);
        num := NumberOfSmallGroups(n);
        printf "\n *** Starting the %o groups of order %o.\n\n", num, n;
        repeat
            H := Current(P);
            if (not IsAbelian(H)) and (not IsPrimePower(Order(H))) then
                phi, G := MinimalDegreePermutationRepresentation(H);
                G;
                for p in PrimeDivisors(Order(G)) do
                    Syl := SylowSubgroup(G, p);
                    if not IsCyclic(Syl) then
                        if not HasNormalpComplementSpecialCase(G,p) then
                            V := ExtTestAuxiliary(G, p, N);
                            if Type(V) eq List then
                                print("We have found a counter-example: ");
                                return(V);
                            end if;
                        end if;
                    end if;
                end for;
            end if;
            Advance(~P);
        until IsEmpty(P);
    end for;
    return true;
end function;
```

Entering the command ExtTestSimples(n,a,b); in Magma, it checks conjecture 5.1 for all finite groups of order a up to order b over a minimal splitting field of characteristic $p$, when $p$ divides the group order. As an example, entering

```
U := ExtTestSimples(5,30,40);U;
```

verifies that conjecture $S_5$ is true for all groups of order 30 up to group order 40. In the Magma code, we used several theoretical results that were proven in the main text, e.g. that the conjecture is true when the $p$-Sylow subgroup is cyclic (corresponding to the representation-finite case) or for groups with a normal $p$-complement.

## References

[AM] Adem, A.; Milgram, R. J.: *Cohomology of finite groups.* Grundlehren der Mathematischen Wissenschaften 309, Springer (2004).

[AnFul] Anderson, F.; Fuller, K.: *Rings and Categories of Modules.* Graduate Texts in Mathematics, Volume 13, Springer-Verlag, 1992.

[ARS] Auslander, M.; Reiten, I.; Smalo, S.: *Representation Theory of Artin Algebras* Cambridge Studies in Advanced Mathematics, 36. Cambridge University Press, Cambridge, 1997. xiv+425 pp.

[Ben] Benson, D.J.: *Representations and cohomology. I: Basic representation theory of finite groups and associative algebras.*

[Br] Brown, K.: *Cohomology of groups.* Graduate Texts in Mathematics 87, Springer (1982). Cambridge Studies in Advanced Mathematics, Cambridge University Press (1991).

[BCP] Bosma,W.; Cannon, J.; Playoust, C.:*The Magma algebra system. I. The user language.*, J. Symbolic Comput., 24 (1997), 235-265.

[CTVZ] Carlson, J.; Townsley, L.; Valeri-Elizondo, L.; Zhang, M.: *Cohomology rings of finite groups. With an appendix: Calculations of cohomology rings of groups of order dividing 64 by Carlson, Valeri-Elizondo and Zhang.* Algebra and Applications, 3. Kluwer Academic Publishers, Dordrecht, 2003. xvi+776 pp.

[CE] Cartan, H.; Eilenberg, S.: *Homological algebra.* Princeton Landmarks in Mathematics, Princeton University Press (1999).

[E] Erdmann, K.: *Blocks of tame representation type and related algebras*, Springer Lecture Notes in Mathematics 1428, 1990.

[ES] Erdmann, K.; Skowroński, A.: *Periodic algebras.* Trends in representation theory of algebras and related topics. Proceedings of the 12th international conference on representations of algebras and workshop (ICRA XII) 2007. European Mathematical Society, EMS Series of Congress Reports, 201-251 (2008).

[ES2] Erdmann, K.; Skowroński, A.: *Hybrid algebras.* arxiv 2103.05963.

[ES3] Erdmann, K.; Skowroński, A.: *The periodicity conjecture for blocks of group algebras.* Colloq. Math. 138, No. 2, 283-294 , 2015.

[ES4] Erdmann, K.; Skowroński, A.: *The stable Calabi-Yau dimension of tame symmetric algebras.* J. Math. Soc. Japan 58 (2006), no. 1, 97-128.

[ES5] K. Erdmann, A. Skowroński, *Weighted surface algebras: general version.* J. Algebra 544 (2020), 170-227.

[GR] Gabriel, P.; Riedtmann, C.: *Group representations without groups.* Commentarii mathematici Helvetici 1979, Volume: 54, page 240-287.

[GKM] Geranios, H.; Kleshchev, A.; Morotti, L.: *On self-extensions of irreducible modules over symmetric groups.* Trans. Am. Math. Soc. 375, No. 4, 2627-2676 (2022).

[G] Green, J. A.: *Walking around the Brauer tree.* Journal of the Australian Mathematical Society , Volume 17 , Issue 2 , 1974 , 197-213.

[Gus] Gustafson, W.: *Global dimension in serial rings.* Journal of Algebra, Volume 97, pages 14-16, 1985.

[H] Holm, T.: *Derived equivalence classification of algebras of dihedral, semidihedral, and quaternion type.* J. Algebra 211 (1999), no. 1, 159-205.

[HB] Huppert, B.; Blackburn, N.: *Finite groups. II.* 44. Grundlehren der Mathematischen Wissenschaften, 242. Springer-Verlag, Berlin-New York, 1982.

[ILP] Igusa, K.; Liu, S.; Paquette, C.: *A proof of the strong no loop conjecture.* Advances in Mathematics Volume 228, Issue 5, 1 December 2011, Pages 2731-2742.

[LM] Liu, S.; Morin, J.: *The strong no loop conjecture for special biserial algebras.* Proc. Amer. Math. Soc. 132 (2004), no. 12, 3513-3523.

[Mar] Marczinzik, R.: *Upper bounds for the dominant dimension of Nakayama and related algebras.* Journal of Algebra.

[Mar2] Marczinzik, R.: *On a conjecture about dominant dimensions of algebras.* Journal of Algebra Volume 565, 2021, Pages 582-597.

[QPA] The QPA-team, QPA - Quivers, path algebras and representations - a GAP package, Version 1.25; 2016 (https://folk.ntnu.no/oyvinso/QPA/).

[Sko] Skowroński, A.: *Minimal representation-infinite artin algebras.* Mathematical Proceedings of the Cambridge Philosophical Society , Volume 116 , Issue 2 , 1994 , 229-243.

[SkoYam] Skowroński, A.; Yamagata, K.: *Frobenius Algebras I: Basic Representation Theory.* EMS Textbooks in Mathematics, 2011.

[S] Swan, R.G.: *Periodic resolutions for finite groups.* Ann. of Math. 72 (1960), 267-291.

[We] Webb, P.: *A course in finite group representation theory.* Cambridge Studies in Advanced Mathematics 161. Cambridge University Press 2016.

(Böhmler) Leibniz University Hannover, Institute of Algebra, Number Theory and Discrete Mathematics, Welfengarten 1, 30167 Hannover, Germany
*E-mail address*: `boehmler@math.uni-hannover.de`

(Erdmann) Mathematical Institute, Oxford University, ROQ, Oxford OX2 6GG, United Kingdom
*E-mail address*: `erdmann@maths.ox.ac.uk`

(Klász) Mathematical Institute of the University of Bonn, Endenicher Allee 60, 53115 Bonn, Germany
*E-mail address*: `klasz@math.uni-bonn.de`

(Marczinzik) Mathematical Institute of the University of Bonn, Endenicher Allee 60, 53115 Bonn, Germany
*E-mail address*: `marczire@math.uni-bonn.de`